\documentclass[oupthm]{CUP-JNL-BCM}%

\usepackage{graphicx}
\usepackage{multicol}
\usepackage{amsmath,amssymb,amsfonts, mathtools}
\usepackage{mathrsfs}
\usepackage{amsthm}
\usepackage{rotating}
\usepackage{appendix}
\usepackage[numbers]{natbib}
\usepackage{changepage,afterpage}
\usepackage{etruscan}
\usepackage{tikz, siunitx, newunicodechar}
\usepackage{booktabs,multirow}
\usepackage{tcolorbox}
\usepackage{bm,bbm}%
\usepackage{braket}

\usepackage[utf8]{inputenc}

\usepackage{array}
\usepackage{xcolor}
\usepackage{colortbl}

\allowdisplaybreaks[4]

\theoremstyle{oupplain}
\newtheorem{theorem}{Theorem}[section]

\newtheorem{corollary}[theorem]{Corollary}
\theoremstyle{oupdefinition}

\theoremstyle{oupremark}

\newtheorem{example}[theorem]{Example}
\theoremstyle{oupproof}

\numberwithin{equation}{section}

\articletype{RESEARCH ARTICLE}

\begin{document}
\setlength{\abovedisplayskip}{3pt}
\setlength{\belowdisplayskip}{3pt}
\renewcommand{\today}{}
\renewcommand{\currenttime}{}

\begin{Frontmatter}
\title{Majorization between symplectic spectra of positive semidefinite matrices}
\author{Temjensangba}
\address{
\orgname{Department of Mathematics and Computing, Indian Institute of Technology (ISM) Dhanbad},
 \orgaddress{\street{Dhanbad, Jharkhand}, 
 \state{India},
  \postcode{826 004}}
  \email{temjensangba111@gmail.com}}
\author{Hemant K. Mishra,}
\address{
\orgname{Department of Mathematics and Computing, Indian Institute of Technology (ISM) Dhanbad}, 
\orgaddress{\street{Dhanbad, Jharkhand}, \state{India}, \postcode{826 004}}
\email{hemantmishra1124@gmail.com}}
\author{Niloy Paul}
\address{
\orgname{Department of Mathematics and Computing, Indian Institute of Technology (ISM) Dhanbad}, 
\orgaddress{\street{Dhanbad, Jharkhand}, \state{India}, \postcode{826 004}}
\email{niloy.paul009@gmail.com}}
\abstract{
    Given $2n \times 2n$ real symmetric positive semidefinite matrix $A$ with symplectic kernel, there exists a real $2n \times 2n$ \emph{symplectic matrix} $M$ such that $M^TAM= D \oplus D$, where $D$ is an $n \times n$ non-negative diagonal matrix which is unique up to permutation of its diagonal entries.
    The diagonal entries of $D$ are called the \emph{symplectic eigenvalues} or symplectic spectrum of $A$. 
    In this work, we investigate some majorization and weak supermajorization relations between the symplectic spectra of two positive semidefinite matrices. 
    More explicitly, suppose $A$ and $B$ are $2n \times 2n$ real symmetric positive semidefinite matrices with symplectic kernels.
    We show that if the symplectic spectrum of $A$ is majorized by the symplectic spectrum of $B$, then $A$ lies in the convex hull of the symplectic orbit of $B$. 
    We also establish that only a weak converse of this statement holds; i.e., if $A$ lies in the convex hull of the symplectic orbit of $B$ then the symplectic spectrum of $A$ is \emph{weakly supermajorized} by the symplectic spectrum of $B$. 
    Several consequences of our results are also presented. 
    Our methods make use of well-known connections between the theory of majorization, doubly stochastic, doubly superstochastic, and symplectic matrices.}
\end{Frontmatter}
\section{Introduction}
    A significant consequence of the $1936$ seminal work of J.~Williamson~\cite{williamson1936algebraic} is that any $2n \times 2n$ real symmetric positive definite matrix can be diagonalized by means of a symplectic matrix. 
   The diagonal entries of the associated diagonal matrix are known as the symplectic eigenvalues or the symplectic spectrum of the positive definite matrix. 
   Generally known as Williamson's theorem, the result has been generalized to positive semidefinite matrices with symplectic kernels \cite{jm, son2022symplectic}. 
    Symplectic spectrum has become an object of interest in various fields such as quantum mechanics \cite{de2006symplectic}, quantum information theory \cite{serafini2023quantum}, and Hamiltonian dynamics \cite{hofer2012symplectic}.
    In particular, Williamson's theorem plays a privileged role in the rich mathematical formalism of bosonic Gaussian information theory \cite{weedbrook2012gaussian}.
    Advancements in the theory of symplectic spectrum holds great potential to expand our understanding of bosonic Gaussian quantum states.
    A significant progress in the symplectic spectral theory parallel to the classic spectral theory has been made in the past few decades.
    This includes symplectic analogs of classic matrix analysis results such as Weyl's ineqaulities \cite{bhatia2021variational}, Lidskii's theorem \cite{jm}, Schur--Horn theorem \cite{bhatia2020schur, huang2023new}, Horn's conjecture \cite{paradan2022horn}, and differentiability and analyticity properties \cite{jm, mishra2020first}.
    Recognizing the fruitfulness and utility of developing symplectic analogs of classic eigenvalue results, our current work investigates symplectic counterparts of a well-known result from the majorization theory of eigenvalues. We deliberate it in the following paragarph.
 
    Given two $n \times n$ Hermitian matrices $A$ and $B$, the spectrum of $A$ is majorized by the spectrum of $B$ if and only if $A$ lies in the convex hull of the unitary orbit of $B$, i.e., there exist $n \times n$ unitary matrices $U_1, \ldots, U_k$ and positive numbers $p_1, \ldots, p_k$ satisfying ${\sum_{j=1}^k p_j=1}$ such that
    \begin{equation}\label{eq:majorizationofhermitianmatrices}
        A = \sum_{j = 1}^k p_j U^*_j B U_j,
    \end{equation}
    where $U^*_j$ denotes the complex conjugate transpose of $U_j$ \cite[Theorem~2.2]{alberti1982stochasticity}.  In this work, we establish the following symplectic analog of the result stated above. 
    Suppose $A$ and $B$ are $2n \times 2n$ real symmetric positive semidefinite matrices with symplectic kernels.
    We show in Theorem~\ref{thm:maj_implies_convex_symplectic} that if the symplectic spectrum of $A$ is majorized by the symplectic spectrum of $B$, then $A$ lies in the convex hull of the \emph{symplectic orbit} of $B$, i.e.,  there exist $2n \times 2n$ real symplectic matrices $M_1, \ldots, M_k$ and positive numbers $p_1, \ldots, p_k$ satisfying ${\sum_{j=1}^k p_j=1}$ such that
    \begin{equation} \label{eq:desirable}
        A = \sum_{j=1}^k p_j M_j^T B M_j.
    \end{equation}
    Example~\ref{firstexample} illustrates that the converse of the above statement does not hold in general.
    Nevertheless, we show in Theorem~\ref{thm:theorem2} that a weak converse holds; i.e., if $A$ lies in the convex hull of the symplectic orbit of $B$ then the symplectic spectrum of $A$ is \emph{weakly supermajorized} by the symplectic spectrum of $B$. 
    We also report several consequences of our main results which are interesting in their own right.
    
    The paper is organized as follows. In Section~\ref{sec:preliminaries}, we provide a brief background on the topics needed during the course of the paper. 
    We state the main results, along with interesting consequences, in Section~\ref{sec:mainresults}. 
    The proofs are provided in Section~\ref{sec:proofs}.
    
\section{Preliminaries}\label{sec:preliminaries}
    We begin by setting some notations.
    Let $\operatorname{S}_n$ denote the permutation group on $\{1,\ldots, n\}$. 
    For every $\pi \in \operatorname{S}_n$, we denote by $P_\pi$ the $n \times n$ permutation matrix whose $(i, j)$th entry is $\delta_{\pi(i), j}$.
    We will denote by $\operatorname{Pd}(2n)$ and $\operatorname{Psd}(2n)$ the sets of all $2n \times 2n$ real symmetric positive definite and semidefinite matrices, respectively.

    Consider the matrix $J \coloneqq \begin{bsmallmatrix} 0 & I_n \\ -I_n & 0 \end{bsmallmatrix}$, where $I_n$ is the identity matrix of size $n$.
    A ${2n \times 2n}$ real matrix $M$ is said to be symplectic if it satisfies $M^TJM = J$. 
    The set of all $2n \times 2n$ real symplectic matrices, denoted by $\operatorname{Sp}(2n)$, forms a group under matrix multiplication and is known as the symplectic group \cite{dms}.
    A linear subspace $\mathcal{W}$ of $\mathbb{R}^{2n}$ is said to be a symplectic subspace if for every $0 \neq u \in \mathcal{W}$ there exists $v \in \mathcal{W}$ such that $u^T Jv \neq 0$.
    A generalization of Williamson's theorem \cite{williamson1936algebraic} states that for any $A \in \operatorname{Psd}(2n)$ with symplectic kernel, there exists $M \in \operatorname{Sp}(2n)$ such that
    \begin{align}\label{eq:Williamson_decomposition}
        M^T A M 
            &=
            \begin{bmatrix}
             D & 0 \\
             0 & D
            \end{bmatrix},
    \end{align}
    where $D$ is an $n \times n$ diagonal matrix with non-negative diagonal entries.
    The diagonal entries of $D$ are unique up to permutation and are called the symplectic eigenvalues of $A$ \cite{jm, son2022symplectic}. In addition, if the symplectic matrix $M$ is orthogonal in \eqref{eq:Williamson_decomposition}, $A$ is said  to be orthosymplectically diagonalizable in the sense of Williamson's theorem. This condition is equivalent to $AJ = JA$ \cite{son2021computing, kamat2024simultaneous}.
    
     Let $A \in \operatorname{Psd}(2n)$ with symplectic kernel.
     We denote by $\operatorname{Sp}(2n, A)$ the set of $2n \times 2n$ real symplectic matrices that diagonalize $A$ in the sense of \eqref{eq:Williamson_decomposition}.
     The symplectic orbit of $A$ is the set $\{ M^T A M : M \in \operatorname{Sp}(2n) \}$.
    We denote by $d(A)$ the $n$-vector with entries given by the symplectic eigenvalues ${d_1(A) \leq \cdots \leq d_n(A)}$ of $A$ arranged in the non-decreasing order.
    Let $M \in \operatorname{Sp}(2n, A)$, and let $x_1,\ldots, x_n, y_1,\ldots, y_n$ be the columns of $M$.
    The following representation of $A$ will play a key role in the proof of Theorem~\ref{thm:maj_implies_convex_symplectic}, which can be verified using the definition of $M$ and the diagonalization \eqref{eq:Williamson_decomposition}:
        \begin{align} 
            A&= \sum_{i=1}^n d_i(A)J\left[ x_i x_i^T +y_i y_i^T \right]J^T. \label{eq:anotherWilliamsonsdecompositionofA}
        \end{align}

\subsection{Relevant concepts from majorization theory} \label{majorization}
    For any $x =(x_1,\ldots, x_n) \in \mathbb{R}^n$, we denote by $x^{\downarrow}=(x_1^{\downarrow},\ldots, x_n^\downarrow)$ the vector obtained by rearranging the entries of $x$ in the non-decreasing order.
    We similarly define $x^\uparrow$.
    Let $x, y \in \mathbb{R}^n$. 
    We say that $x$ is weakly supermajorised by $y$, in symbols $x \prec^w y$, if
    \begin{equation} \label{eq:weaksupermajorization}
    \sum_{j = 1}^k x_j^\uparrow \geq \sum_{j=1}^k y_j^\uparrow \qquad \text{ for } 1 \leq k \leq n.
    \end{equation}
     If the equality holds in \eqref{eq:weaksupermajorization} for $k =n$, we say that $x$ is majorized by $y$ and written in symbols as $x \prec y$.
     An $n \times n$ real matrix with non-negative entries is said to be doubly stochastic if each row, as well as each column, sums to $1$. An $n \times n$ real matrix $S = [s_{ij}]$ is said to be doubly superstochastic if there exists some doubly stochastic matrix $E = [e_{ij}]$ such that $s_{ij} \geq e_{ij}$ for all $1 \leq i, j \leq n$. 
    We have $x \prec y$ if and only if $x = Ey$ for some $n \times n$ doubly stochastic matrix $E$ \cite[Theorem~2.B.2]{marshall1979inequalities}.
    Also, if $x, y$ are non-zero vectors with only non-negative entries then $x \prec^w y$ if and only if $x = Sy$ for some $n \times n$ doubly superstochastic matrix $S$ \cite[Proposition~1.A.5]{marshall1979inequalities}. 
    Additionally, we write $x \prec_+^w y$ if all the entries of $S$ are strictly positive.

    We will also use the following connection between symplectic matrices and doubly superstochastic matrices.
    Let $M \in \operatorname{Sp}(2n)$, written in the block form
    \begin{align}
        M & = \begin{bmatrix}
            P & Q \\
            R & S
        \end{bmatrix},
    \end{align}
    where each block is of size $n \times n$.
    Then the matrix 
    \begin{align}\label{eq:dsupstochastic_associated_symplectic}
    \widetilde{M} \coloneqq \frac{1}{2} (P \circ P + Q \circ Q + R \circ R + S \circ S)
    \end{align}
    is doubly superstochastic \cite[Theorem~6]{bhatia2015symplectic}.
    Here $\circ$ denotes the Hadamard product of matrices.

    \subsection{Certain vectors associated with a positive definite matrix}
    \label{symplecticdiagonals}

    Let $A = \begin{bsmallmatrix} A_{11} & A_{12} \\ A_{12}^T & A_{22} \end{bsmallmatrix} \in \operatorname{Pd}(2n)$ with the $n \times n$ blocks $A_{11}, A_{12}, A_{22}$. 
    Denote by $\Delta_{11}, \Delta_{12}$, and $\Delta_{22}$ the $n$-vectors consisting of the diagonal entries of $A_{11}, A_{12}$, and $A_{22}$, respectively. 
    Define
    \begin{align}
    \Delta_c(A) &:= \frac{\Delta_{11} + \Delta_{22}}{2}, \label{eq:average_diagonal} \\
    \Delta_s(A) &:= \sqrt{\Delta_{11} \cdot \Delta_{22}}, \\
    \Delta_h(A) &:= \sqrt{\frac{\Delta_{11}^2 + \Delta_{22}^2}{2}}, \\
    \Delta_w(A) &:= \sqrt{\frac{\Delta_{11}^2 + \Delta_{22}^2 + 2 \Delta_{12}^2}{2}} \label{eq:w_diagonal}. 
    \end{align}
    Here $\Delta_{11} \cdot \Delta_{22}$ means the entry-wise multiplication of the vectors.
    It is known that the following weak supermajorization relation holds:
    \begin{align}
        \Delta_\xi(A) \prec^w d(A)
    \end{align}
    for all $\xi \in \{c, s, h, w\}$.
    See \cite{bhatia2020schur, huang2023new}.
    The following symplectic matrix associated with $A$ will be useful later:
    \begin{equation} \label{eq:Mforc}
        M = \operatorname{diag} \Big(\sqrt[4]{\Delta_{11}^{-1} \cdot \Delta_{22}} \Big) \oplus \operatorname{diag} \Big( \sqrt[4]{\Delta_{22} \cdot \Delta_{11}^{-1}} \Big). 
    \end{equation}
    Again, algebraic operations in $\Delta_{11}^{-1}$ and $\Delta_{11}^{-1} \cdot \Delta_{22}$ are understood in the entry-wise sense.
    
    \subsection{Symplectic direct sum and symplectic pinching} \label{pinching}
    Given any square matrices $X_1, \ldots, X_k$, we denote their direct sum by $\oplus X_i$.
    For any square matrix $X$, let $\mathscr{C}(X)$ denote a pinching of $X$, which is the matrix obtained by retaining some diagonal blocks and making all the other entries zero.
    We now recall symplectic analogs of these notions.
    
    Let $m_1, \ldots, m_k \in \mathbb{N}$ such that $m_1 + \cdots + m_k = n$.
    For each $i \in\{ 1,\ldots k\}$, let $A_i = \begin{bsmallmatrix}
        E_i & F_i \\ G_i & H_i
    \end{bsmallmatrix}$
     be a $2m_i \times 2m_i$ real matrix with $m_i \times m_i$ blocks $E_i, F_i, G_i, H_i$. 
     Recall from \cite{bhatia2015symplectic} that the symplectic direct sum of $A_1, \ldots, A_k$ is defined as
    \begin{equation}
    \oplus^s A_i := \begin{bmatrix}
        \oplus E_i & \oplus F_i \\ \oplus G_i & \oplus H_i 
    \end{bmatrix}. \label{eq:sdirectsum}
    \end{equation}
       
   Let $A = \begin{bsmallmatrix} E & F \\ G & H  \end{bsmallmatrix}$ be a real matrix of size $2n \times 2n$ with $n \times n$ blocks $E,F,G,$ and $H$. 
   The symplectic pinching or $s$-pinching of $A$ is defined as
    \begin{equation}
    \mathscr{C}^s(A) := \begin{bmatrix} \mathscr{C}(E) & \mathscr{C}(F) \\
    \mathscr{C}(G) & \mathscr{C}(H) \end{bmatrix}. \label{eq:spinching}
    \end{equation}
    Further, suppose $\mathscr{C}(E) = \oplus E_i, \mathscr{C}(F) = \oplus F_i, \mathscr{C}(G) = \oplus G_i$, and $\mathscr{C}(H) = \oplus H_i$.
    Then we have
    \begin{equation} 
        \mathscr{C}^s(A) = \oplus^s A_i, \label{eq:spinching=sdirectsum}
    \end{equation}
    where $A_i \coloneqq \begin{bsmallmatrix} 
            E_i & F_i \\ G_i & H_i 
        \end{bsmallmatrix}$ for all $1 \leq i \leq k$.

\section{Results} \label{sec:mainresults}
    We begin by stating the first main result. 
    \begin{theorem} \label{thm:maj_implies_convex_symplectic}
        Let $A, B \in \operatorname{Psd}(2n)$ with symplectic kernels. 
        If $d(A) \prec d(B)$, then there exist a family of symplectic matrices $M_\pi$ for $\pi \in \operatorname{S}_n$ and a probability vector $( p(\pi) )_{\pi \in S_n}$ such that 
        \begin{equation} \label{eq:equation_in_first_result} 
            A = \sum_{\pi \in S_n} p(\pi) M_\pi^T B M_\pi. 
        \end{equation}
    \end{theorem}

        We note in passing a modified version of Theorem~\ref{thm:maj_implies_convex_symplectic} obtained by replacing the hypothesis $d(A) \prec d(B)$ by $d(A) \prec_+^w d(B)$.
    \begin{theorem} \label{thm:modifiedversion}
        Let $A, B \in \operatorname{Psd}(2n)$ such that their kernels are symplectic subspaces of $\mathbb{R}^{2n}$. If $d(A) \prec_+^w d(B)$, then there exists a probability vector $(p(\pi))_{\pi \in \operatorname{S}_n}$ such that
        \begin{equation} \label{eq:modifiedequation}
            A = \sum_{\pi \in \operatorname{S}_n} p(\pi) \Big( K_\pi D_\pi M \Big)^T B \Big(K_\pi D_\pi M \Big),
        \end{equation}
        where $K_\pi$ and $M$ are symplectic matrices and each $D_\pi$ is a positive diagonal matrix.
    \end{theorem}
    The following example is a striking illustration of the fact that the converse of Theorem~\ref{thm:maj_implies_convex_symplectic} does not hold in general.

\begin{example} \label{firstexample}
Take $B = I_4$, and the symplectic matrices $M_1 = I_4$ and 
$M_2 = 2 I_2 \oplus \frac{1}{2} I_2$.
Consider
\begin{equation}
A = \frac{1}{2} M_1^TBM_1 + \frac{1}{2} M_2^TBM_2 = \dfrac{5}{2} I_2 \oplus \frac{5}{8} I_2.
\end{equation}
We have $d(A) = \left(\frac{5}{4}, \frac{5}{4} \right)$ and $d(B) = \left(1, 1\right)$.
The weak supermajorization $d(A) \prec^w d(B)$ holds but the majorization $d(A) \prec d(B)$ does not hold.
\end{example} 

As the second main result, we show in the next theorem that a ``weak'' converse of Theorem~\ref{thm:maj_implies_convex_symplectic} is true.
\begin{theorem} \label{thm:theorem2}
Let $A, B \in \operatorname{Psd}(2n)$ such that their kernels are symplectic subspaces of $\mathbb{R}^{2n}$. 
If $A$ lies in the convex hull of the symplectic orbit of $B$, then $d(A) \prec^w d(B)$.
\end{theorem}

    We report some interesting consequences of our main results for positive definite matrices in the following corollaries.
    
    \begin{corollary} \label{thm:firstcorollary}
        Let $A \in \operatorname{Pd}(2n)$ such that $AJ = JA$, and fix $\xi \in \{ c, w, h \}$. 
        Every $B \in \operatorname{Pd}(n)$ satisfying $d(B) = \Delta_\xi(A)^\uparrow$ lies in the convex hull of the symplectic orbit of $A$.
        In particular, $\operatorname{diag} \left( \Delta_\xi (A) \right) \oplus \operatorname{diag} \left( \Delta_\xi (A) \right)$ lies in the convex hull of the symplectic orbit of $A$.
    \end{corollary}

    \begin{corollary} \label{thm:secondcorollary}
     Let $A \in \operatorname{Pd}(2n)$ such that $(MAM)J = J(MAM)$, where $M$ is as defined in \eqref{eq:Mforc}. Every ${B \in \operatorname{Pd}(2n)}$ satisfying $d(B) = \Delta_s(A)^\uparrow$ lies in the convex hull of the symplectic orbit of $A$.
    In particular, $\operatorname{diag} \left( \Delta_{11} \right) \oplus \operatorname{diag} \left( \Delta_{22} \right)$ lies in the convex hull of the symplectic orbit of $A$, where $\Delta_{11}$ and $\Delta_{22}$ are vectors associated with $A$ as defined in Section~\ref{symplecticdiagonals}.
    \end{corollary}

    \begin{corollary} \label{thm:thirdcorollary}
        For $A \in \operatornamewithlimits{Pd}(2n)$, let $\mathscr{C}^s(A) = \oplus^s A_i$ as defined in \eqref{eq:spinching=sdirectsum}. 
        Choose and fix ${M_i \in \operatorname{Sp}(2m_i, A_i)}$, and let $M \coloneqq \oplus^s M_i$. 
        If $(M^TAM)J = J(M^TAM)$, then $\mathscr{C}^s(A)$ lies in the convex hull of the symplectic orbit of $A$.
    \end{corollary}

\section{Proofs} \label{sec:proofs}
    Our proofs of Theorem \ref{thm:maj_implies_convex_symplectic} and Theorem \ref{thm:theorem2} are inspired from the proof of the classic counterpart result of eigenvalues presented in J. Watrous's treatise \cite[Theorem~4.33]{watrous2018theory}.
     
    To avoid trivial cases, we will assume that the positive semidefinite matrices $A$ and $B$ considered in the statements of the theorems are all non-zero.
 
\subsection{Proof of Theorem \ref{thm:maj_implies_convex_symplectic}.}
     We know by \cite[Theorem~1.3]{ando1989} that the hypothesis $d(A) \prec d(B)$ implies that there exists a doubly stochastic matrix $E$ such that
        \begin{equation}
            d(A) = Ed(B).
        \end{equation}
        By Birkhoff-von Neumann theorem \cite{birkhoff1946three}, there exists a probability vector $\left(p(\pi) \right)_{\pi \in \operatorname{S}_n}$ such that
        \begin{equation}
            E = \sum_{\pi \in \operatorname{S}_n} p(\pi)P_{\pi}.
        \end{equation}
        Therefore
        \begin{equation}
          d(A) = \sum_{\pi \in \operatorname{S}_n} p(\pi)P_{\pi}d(B), 
        \end{equation}
        which implies that for all $1 \leq i \leq n$,
        \begin{equation}
            d_i(A) =\sum_{\pi \in \operatorname{S}_n} p(\pi)d_{\pi(i)}(B).
        \end{equation}
        Let $[x_1, \ldots, x_n, y_1, \ldots, y_n] \in\operatorname{Sp}(2n, A)$ and $[u_1, \ldots, u_n, v_1, \ldots, v_n] \in\operatorname{Sp}(2n, B)$ be arbitrary.
        It follows from \eqref{eq:anotherWilliamsonsdecompositionofA} that
        \begin{align}
            A&=  \sum_{i=1}^n \sum_{\pi \in \operatorname{S}_n} p(\pi)d_{\pi(i)}(B)J\left[ x_i x_i^T +y_i y_i^T \right]J^T.
        \end{align}
        For every $\pi \in \operatorname{S}_n$, choose
        \begin{align}
           N_{\pi}&=\sum_{i=1}^n  \left[  u_{\pi(i)} x_i^T+ v_{\pi(i)}y_i^T\right].
        \end{align}
        Let us now consider 
        \begin{align}
            \sum_{\pi \in \operatorname{S}_n} p(\pi) N_{\pi}^T B N_{\pi}
                &= \sum_{\pi \in \operatorname{S}_n} p(\pi) N_{\pi}^T \sum_{i=1}^n \left[Bu_{\pi(i)}x_i^T + Bv_{\pi(i)}y_i^T \right] \label{eq:similarsteps} \\
                &= \sum_{\pi \in \operatorname{S}_n} p(\pi) N_{\pi}^T \sum_{i=1}^n \left[d_{\pi(i)}(B)J\left(v_{\pi(i)}x_i^T - u_{\pi(i)}y_i^T\right) \right] \\
                &= \sum_{\pi \in \operatorname{S}_n} p(\pi)  \sum_{i=1}^n \left[d_{\pi(i)}(B)\left(N_{\pi}^TJv_{\pi(i)}x_i^T - N_{\pi}^TJu_{\pi(i)}y_i^T\right) \right].  \label{eq:intermediate}
        \end{align}
        We observe that $N_{\pi} = [u_{\pi(1)}, \ldots, u_{\pi(n)}, v_{\pi(1)}, \ldots, v_{\pi(n)}] [x_1, \ldots, x_n, y_1, \ldots, y_n]^T$ is a symplectic matrix. Therefore, we obtain
        \begin{equation}
        N_{\pi}^T J v_{\pi(i)} = x_i.
        \end{equation}
        Similarly 
        \begin{equation}
            N_{\pi}^T J u_{\pi(i)} = - y_i.
        \end{equation}
        Consequently, \eqref{eq:intermediate} becomes
        \begin{align}
        \sum_{\pi \in \operatorname{S}_n} p(\pi) N_{\pi}^T B N_{\pi} &= \sum_{\pi \in \operatorname{S}_n} p(\pi)  \sum_{i=1}^n \left[d_{\pi(i)}(B)\left(x_ix_i^T + y_iy_i^T\right) \right] \\
                &= \sum_{i=1}^n \sum_{\pi \in \operatorname{S}_n} p(\pi) d_{\pi(i)}(B) \left(x_{i}x_i^T + y_{i}y_i^T\right) \\
                &= \sum_{i=1}^n d_i(A)\left(x_{i}x_i^T + y_{i}y_i^T\right) \\
                &= J^T A J.
        \end{align}
        We thus get
        \begin{align} \label{eq:A_in_terms_of_B}
            A &= \sum_{\pi\in \operatorname{S}_n} p(\pi) (N_{\pi}J^T)^T B (N_{\pi} J^T).
        \end{align}
        By taking $M_\pi = N_{\pi}J^T$, we obtain the desired expression \eqref{eq:equation_in_first_result}.
        \qed

\subsection{Proof of Theorem \ref{thm:modifiedversion}.}
    The proof of this theorem is a modified version of the proof of the previous Theorem~\ref{thm:maj_implies_convex_symplectic} obtained by invoking Sinkhorn's~\cite{sinkhorn1964relationship} result.
 
    By the assumption $d(A) \prec_+^w d(B)$, there exists some doubly superstochastic matrix $K$ having all entries positive such that $d(A) = Kd(B)$. Since all the entries in $K$ are positive, there exist $n \times n$ positive diagonal matrices $D_1, D_2$ and a doubly stochastic matrix $E$ such that $K = D_1 E D_2$ \cite[Theorem~1]{sinkhorn1964relationship}. 
    We know by the Birkhoff-von Neumann theorem that there exists a probability vector $(p(\pi))_{\pi \in \operatorname{S}_n}$ such that
    \begin{equation}
        E = \sum_{\pi \in \operatorname{S}_n} p(\pi) P_\pi.
    \end{equation}
    Therefore, $K$ can be written as
    \begin{equation}
    K = \sum_{\pi \in \operatorname{S}_n} p(\pi) D_1 P_\pi D_2.
    \end{equation}
    Suppose $D_1 = \operatorname{diag}(a_1, \ldots, a_n)$ and $D_2 = \operatorname{diag}(b_1, \ldots, b_n)$. For all $1 \leq i \leq n$, it can be verified that
    \begin{equation}
        d_i(A) = \sum_{\pi \in \operatorname{S}_n} p(\pi) a_i b_{\pi(i)} d_{\pi(i)} (B).
    \end{equation}
    Let $[x_1, \ldots, x_n, y_1, \ldots, y_n] \in\operatorname{Sp}(2n, A)$ and $[u_1, \ldots, u_n, v_1, \ldots, v_n] \in\operatorname{Sp}(2n, B)$ be arbitrary.
    For every $\pi \in S_n$, choose
    \begin{equation}
      N_\pi = \sum_{i =1}^n \left[ \sqrt{b_{\pi(i)}} u_{\pi(i)} \left( \sqrt{a_i} x_i \right)^T + \sqrt{b_{\pi(i)}} v_{\pi(i)} \left( \sqrt{a_i} y_i \right)^T \right].   
    \end{equation}
    By following the same steps as in the proof of Theorem~\ref{thm:maj_implies_convex_symplectic} from \eqref{eq:similarsteps} onward, we obtain
    \begin{equation}
    A = \sum_{\pi \in \operatorname{S}_n} p(\pi) (K_\pi D_\pi M)^T B (K_\pi D_\pi M),
    \end{equation}
    where
    \begin{align}
        K_\pi &= [u_{\pi(1)}, \ldots, u_{\pi(n)}, v_{\pi(1)}, \ldots, v_{\pi(n)}], \\
        D_\pi &= \operatorname{diag} \left(\sqrt{a_1 b_{\pi(1)}}, \ldots, \sqrt{a_n b_{\pi(n)}}, \sqrt{a_1 b_{\pi(1)}}, \ldots, \sqrt{a_n b_{\pi(n)}} \right). \\
        M &= \left[ x_1, \ldots, x_n, y_1, \ldots, y_n \right] J^T.
    \end{align}
    \qed

\subsection{Proof of Theorem \ref{thm:theorem2}.}
    Let $[x_1, \ldots, x_n, y_1, \ldots, y_n] \in\operatorname{Sp}(2n, A)$ and $[u_1, \ldots, u_n, v_1, \ldots, v_n] \in\operatorname{Sp}(2n, B)$ be arbitrary.
    We have for $1 \leq j\leq n$ that
        \begin{equation}\label{eq:symp_eig_A_alternate_exp}
            d_j(A) 
                = \frac{1}{2} \left( x_j^T A x_j + y_j^T A y_j \right).
        \end{equation}
    Assume that
        \begin{equation}\label{eq:assumption_A_convexhull_B}
            A = \sum_{i=1}^m \lambda_i M_i^T B M_i,
        \end{equation}
    where $\left(\lambda_1,\ldots, \lambda_m\right)$ is a probability vector with strictly positive entries, and $M_1,\ldots, M_m$ are symplectic matrices.
    By substituting \eqref{eq:assumption_A_convexhull_B} into \eqref{eq:symp_eig_A_alternate_exp}, we get
        \begin{align}
            d_j(A) &= \frac{1}{2} \sum_{i=1}^m \lambda_i \Big( x_j^T M_i^T B M_i x_j + y_j^T M_i^T B M_i y_j \Big) \\
            &= \frac{1}{2} \sum_{i=1}^m \lambda_i \sum_{k=1}^n d_k(B) \Big( x_j^T M_i^T J (u_k u_k^T + v_k v_k^T) J^T M_i x_j +  y_j^T M_i^T J (u_k u_k^T + v_k v_k^T) J^T M_i y_j \Big) \\
            &= \frac{1}{2} \sum_{i=1}^m \lambda_i \sum_{k=1}^n d_k(B) \Big( (x_j^T M_i^T J u_k)^2 + (x_j^T M_i^T J v_k)^2 + (y_j^T M_i^T J u_k)^2 + (y_j^T M_i^T J v_k)^2 \Big) \\
            &= \sum_{i=1}^m \lambda_i \sum_{k=1}^n \frac{1}{2} \Big( (x_j^T M_i^T J u_k)^2 + (x_j^T M_i^T J v_k)^2 + (y_j^T M_i^T J u_k)^2 + (y_j^T M_i^T J v_k)^2 \Big) d_k(B). 
        \end{align}
        This gives
        \begin{equation}
        d(A) =\sum_{i=1}^m  \lambda_i \widetilde{S}_i d(B),  
        \end{equation}
        where each $\widetilde{S}_i$ is the doubly superstochastic matrix associated with the symplectic matrix
        \begin{align}
            S_i &= \left[x_1,\ldots, x_n, y_1,\ldots, y_n \right]^T M_i^T J \left[u_1,\ldots, u_n, v_1,\ldots, v_n \right],
        \end{align}
        as defined by \eqref{eq:dsupstochastic_associated_symplectic}.
        Choose $\displaystyle{K = \sum_{i=1}^m  \lambda_i \widetilde{S}_i}$ so that $K$ is a doubly superstochastic matrix satisfying the equality $d(A)= Kd(B)$.
        By Proposition~1.A.5 of \cite{marshall1979inequalities}, we thus obtain $d(A) \prec^w d(B)$.
    \qed

\subsection{Proof of Corollary~\ref{thm:firstcorollary}.}
    We know from \cite[Theorem~5.4]{mishraequality} and \cite[Theorem~1.2]{huangandmishra} that $AJ = JA$ implies $\Delta_\xi (A) \prec d(A)$ for all $\xi \in \{c, w, h \}$.
    So, if $d(B) = \Delta_\xi(A)^\uparrow$ then we have $d(B) \prec d(A)$.
    It thus follows from Theorem~\ref{thm:maj_implies_convex_symplectic} that $B$ lies in the convex hull of the symplectic orbit of $A$.

    In particular, the choice $B= \operatorname{diag} \left( \Delta_\xi (A) \right) \oplus \operatorname{diag} \left( \Delta_\xi (A) \right)$ gives $d(B) = \Delta_\xi(A)^\uparrow$, and hence $B$ lies in the convex hull of the symplectic orbit of $A$.
    \qed

\subsection{Proof of Corollary~\ref{thm:secondcorollary}.}
    The proof follows similar arguments as in the proof of Corollary~\ref{thm:firstcorollary}, where we replace $A$ with $MAM$, and use Theorem~1.1 \cite{huangandmishra}.
    \qed

\subsection{Proof of Corollary~\ref{thm:thirdcorollary}.}
    By Theorem~1.4 of \cite{huangandmishra} we know that $\left( M^TAM \right)J = J \left( M^TAM \right)$ implies $d \left( \mathscr{C}^s(A) \right) \prec d(A)$.
    It then directly follows by Theorem~\ref{thm:maj_implies_convex_symplectic} that $\mathscr{C}^s(A)$ lies in the convex hull of the symplectic orbit of $A$.
    \qed

\section*{Acknowledgements}
    Temjensangba thanks Nagaland University for granting study leave with pay. 
    Hemant K. Mishra acknowledges support from FRS Project No.~MISC~0147. 

\begin{Backmatter}
\bibliographystyle{unsrt}
\bibliography{reference}
\printaddress
\end{Backmatter}

\end{document}